\documentclass[preprint,12pt]{elsarticle}

\usepackage{amsmath, amssymb, amsfonts, amsbsy, latexsym}
\usepackage{graphicx}
\usepackage[all]{xy}
\newcommand{\ud}{\mathrm{d}}

\usepackage{amssymb}
 \usepackage{amsthm}

\journal{J. of Mathematical Analysis and Applications}
\newtheorem{theorem}{Theorem}[section]
\newtheorem{lemma}[theorem]{Lemma}
\newtheorem{corollary}[theorem]{Corollary}

\theoremstyle{definition}
\newtheorem{definition}[theorem]{Definition}

\theoremstyle{remark}

\begin{document}
\allowdisplaybreaks

\begin{frontmatter}


 \title{Mixed Modulation Spaces and Their Application to Pseudodifferential Operators}
 \author[sb]{Shannon~Bishop\corref{cor1}\fnref{label2}}
\ead{sbishop@math.gatech.edu}
 \fntext[label2]{Partially supported by NSF Grant DMS-0806532.}
 \cortext[cor1]{Corresponding author.}
 \address[sb]{School of Mathematics,
Georgia Institute of Technology,
Atlanta, GA 30332 USA}

\begin{abstract}
This paper uses frame techniques to characterize the Schatten class properties of integral operators.  The main result shows that if the coefficients \( \left\{ \langle k, \Phi_{m,n} \rangle    \right\} \) of certain frame expansions of the kernel \(k \) of an integral operator are in \( \ell^{2,p} \), then the operator is Schatten \( p \)-class.  As a corollary, we conclude that if the kernel or Kohn-Nirenberg symbol of a pseudodifferential operator lies in a particular mixed modulation space, then the operator is Schatten \( p \)-class.  Our corollary improves existing Schatten class results for pseudodifferential operators and the corollary is sharp in the sense that larger mixed modulation spaces yield operators that are not Schatten class.

\end{abstract}

\begin{keyword}
frames \sep Gabor transform \sep integral operator \sep Modulation space \sep pseudodifferential operator \sep Schatten class

\MSC[2000] 35S05 \sep 42C15 \sep 47B10
% 35S05 is general theory of pseudodifferential operators
% 42C15 Series of general orthogonal functions, generalized Fourier expansions, nonorthogonal expansions
%47B10 Operators belonging to operator ideals (nuclear, $p$-summing, in the Schatten-von Neumann classes, etc.)

\end{keyword}

\end{frontmatter}

\section{Introduction}
Integral operators arise naturally in many areas of mathematics and science.  Pseudodifferential operators, which are a particular type of integral operator, have appeared widely in the literature of physics, signal processing and differential equations.    An overview of pseudodifferential operators is given in Chapter 14 of \cite{groch}, while more detailed expositions are found in \cite{folland}, \cite{horm2}, and \cite{stein}. Because of the role of pseudodifferential operators in partial differential equations, the smoothness of the Weyl and Kohn-Nirenberg symbols of a pseudodifferential operator has traditionally been used to characterize properties of the operator, with the H\"ormander symbol classes playing key roles.

More recently, pseudodifferential operators have been studied from a time-fre\-quency perspective.  Every pseudodifferential operator is a superposition of time-frequency shifts, and the properties of pseudodifferential operators have been well-described by time-frequency analysis. Results with this flavor appear in  \cite{specasym}, \cite{tach1} and \cite{wong}.  In particular the modulation spaces \(  M^{p,q}_{w}( \mathbb{R}^{d}) \), which are Banach spaces characterized by time-frequency shifts, have been useful symbol spaces for studying continuity and Schatten class properties of pseudodifferential operators.  Using Gabor frames, elements in these spaces can be decomposed into a superposition of time-frequency shifts, and this Gabor frame decomposition of the symbol of a pseudodifferential operator can be used to characterize the properties of the operator.  Results of this type appear in \cite{czaja}, \cite{modbounded}, \cite{heilpdo}, \cite{labate}, \cite{tach2} and \cite{toftpdo}, while modulation spaces appear implicitly in \cite{singval}, \cite{sjo}, \cite{daub}, \cite{horm} and \cite{roch1}.

In this paper we develop a technique for analyzing the kernel of an
integral operator which both generalizes and improves existing time-frequency analysis
techniques of pseudodifferential operators, and in particular yields larger non-smooth classes of Kohn-Nirenberg symbols which ensure that a given pseudodifferential operator is Schatten \(p\)-class.

To obtain our main result, we analyze the slices of the kernel of an integral operator with a
frame.  If these decomposed slices have a certain decay, then the operator is Schatten \( p \)-class.  As a special case, we obtain the following theorem.

\begin{theorem} \label{thm:pclassframe} Suppose \( \left\{ \phi_{m} \right\}_{m \in \Lambda} \) is a frame for \( L^{2}(\mathbb{R}^{d}) \).  Let \( \Phi_{m,n} = \phi_{m} \otimes \overline{\phi_{n}}  \).  If \( A \) is an integral operator with kernel \( k \) and \( p \in [1,2] \) then \( A \) is Schatten \(p\)-class on \( L^{2}(\mathbb{R}^{d}) \) if 
\[ \biggl( \sum_{n \in \Lambda} \biggl( \sum_{m \in \Lambda} \left| \langle k, \Phi_{m,n} \rangle \right|^{2} \biggr)^{\frac{p}{2}}  \biggr)^{\frac{1}{p} } < \infty \] 
\end{theorem}

Analyzing the slices of the kernel as in Theorem \ref{thm:pclassframe} with a Gabor frame in particular gives a time-frequency condition on the kernel which ensures the operator is Schatten \( p\)-class.  We show that this condition holds for kernels belonging to certain Banach spaces \( M(c)_{w}^{p_{1}, p_{2}, \dots, p_{2d}} \) that we call mixed modulation spaces, which are natural generalizations of the traditional modulation spaces \(  M^{p,q}_{w}( \mathbb{R}^{d}) \).  In this paper we show that many of the interesting properties of traditional modulation spaces also hold for mixed modulation spaces.  Furthermore, inclusion of the Kohn-Nirenberg symbol in an appropriate mixed modulation space ensures the corresponding operator is Schatten \( p\)-class.  The relationship between mixed modulation spaces and the kernels and Kohn-Nirenberg symbols of Schatten \( p\)-class operators is summarized in the following theorem.

\begin{theorem} \label{thm:maintf} Let \( A \) be a pseudodifferential operator with kernel \( k \) and Kohn-Nirenberg symbol \( \tau \).  Assume  \( p \in [1,2] \) and set \( 2 = p_{1} = \dots = p_{2d} \) and \( p = p_{2d+1} = \dots = p_{4d} \). For suitable \( c \), if one of \( k, \tau  \) lies in \( M(c)^{p_{1},p_{2},\dots, p_{4d}} \) then so does the other.  In this case \( A \) is Schatten \(p\)-class on \( L^{2}(\mathbb{R}^{d}) \).
\end{theorem}

The strongest known Schatten class result for pseudodifferential operators obtained by time-frequency analysis is Theorem 1.2 of
\cite{modbounded}, which states that if the Weyl symbol or Kohn-Nirenberg symbol of a pseudodifferential operator is in \( M^{2,2}_{v_{s}}(\mathbb{R}^{2d}) \), then the operator is Schatten \( p \)-class if \( p > \frac{2d}{d+s} \) and \( s \geq 0 \).  Although the crux of both Theorem \ref{thm:maintf} and \cite[Theorem 1.2]{modbounded} is time-frequency analysis with Gabor frames, our Theorem  \ref{thm:maintf} is obtained by analyzing the \emph{slices} of the kernel with a Gabor frame, thus permitting a finer control on the properties of the kernel (and, consequently, the symbol).  As a result, we can show that Theorem  \ref{thm:maintf} is stronger than \cite[Theorem 1.2]{modbounded}, in the sense that the mixed modulation space described by Theorem \ref{thm:maintf} strictly contains the space \( M^{2,2}_{v_{s}}(\mathbb{R}^{2d}) \).

The mixed modulation space \( M(c)^{p_{1},p_{2},\dots, p_{4d}} \) is characterized by \( 4d \) decay parameters \(p_{1},p_{2},\dots, p_{4d} \), while the mixed modulation space described by Theorem \ref{thm:maintf} essentially only has two decay parameters.  This disparity suggests that Theorem \ref{thm:maintf} may be extended to a larger mixed modulation space by a more subtle analysis of the kernel of a pseudodifferential operator.  However, this is not the case.  In fact, we show that Theorem \ref{thm:maintf} is sharp in the sense that larger mixed modulation spaces contain kernels and symbols of pseudodifferential operators that are not Schatten \( p \)-class.
%Another approach- also bring how my technique was an improvement of Heil's technique for motivation.

  The paper is organized as follows.  Section 2 contains preliminary and background information.  Section 3 is devoted to the proof of Theorem \ref{thm:pclassframe}.  In Section 4, the definition of mixed modulation spaces \( M(c)_{w}^{p_{1}, p_{2}, \dots, p_{2d}} \) is given and the properties of these spaces are described.   In Section 5, we apply the results of Sections 3 and 4 to pseudodifferential operators.

\section{Preliminaries} 

\subsection{Weight functions} 

\begin{definition}  A locally integrable function \( v: \mathbb{R}^{d} \to [0, \infty) \) is called a \emph{weight function}.  A weight function \( v: \mathbb{R}^{d} \to [0, \infty) \) is \emph{submultiplicative}  if \[ v (z_{1} + z_{2}) \leq v( z_{1}) v( z_{2}) \hspace{2pc} \textrm{ for all } z_{1}, z_{2} \in \mathbb{R}^{d}. \]  A weight function \( v \) has \emph{polynomial growth} if there are \( C, s \geq 0 \) such that \( v(z) \leq C \left( 1 + |z| \right)^{s} \) for all \( z \in \mathbb{R}^{d}\).
\end{definition}

For each \( s \geq 0 \), the function \( v_{s}(z)  = \left( 1 + | z | \right)^{s} \) is a submultiplicative weight function with polynomial growth.

\begin{definition}  Suppose \( w: \mathbb{R}^{d} \to [0, \infty) \) is a weight function and \( v: \mathbb{R}^{d} \to [0, \infty) \) is submultiplicative.   If there is a constant \( C \) such that
\[ w \left(z_{1} + z_{2} \right) \leq C \, v \left( z_{1} \right) w \left( z_{2} \right) \hspace{2pc} \textrm{ for all } z_{1}, z_{2} \in \mathbb{R}^{d}, \] then we call \( w \) a \emph{\( v\)-moderate} weight.
\end{definition}

We will assume throughout this paper that \( v: \mathbb{R}^{d} \to [0, \infty) \) is a submultiplicative weight function of polynomial growth symmetric in each coordinate, i.e.
\( v(x_{1}, \dots, -x_{i}, \dots, x_{d}) = v(x_{1}, \dots, x_{i}, \dots, x_{d}) \) for each \( i = 1,2, \cdots, d \).  We also assume throughout that \( w \) is a \( v\)-moderate weight.

\subsection{Mixed norm spaces}

\begin{definition} Given measure spaces \( \left( X_{i}, \mu_{i} \right) \) and given \( p_{i} \in [1, \infty] \) for \( i = 1, 2, \dots, d \), we let  \( L^{p_{1}, p_{2}, \dots, p_{d}}_{w} \left( X_{1}, X_{2}, \dots, X_{d}, \mu_{1}, \mu_{2}, \dots,  \mu_{d}\right) \) consist of all of the measurable functions \( F: X_{1} \times X_{2} \times \cdots \times X_{d} \to \mathbb{C} \) for which the following norm is finite:
\begin{align*} & \left\| F \right\|_{L^{p_{1}, p_{2}, \dots, p_{d}}_{w}}  \\
& = \biggl( \int_{X_{d}} \dots   \biggl( \int_{X_{1}} \left| F (x_{1}, \dots, x_{d} )   w (x_{1}, \dots, x_{d}) \right|^{p_{1}} \, \ud \mu_{1} (x_{1})  \biggr)^{\frac{p_{2}}{p_{1}}} \dots  \, \ud \mu_{d} (x_{d}) \biggr)^{\frac{1}{p_{d}}},
\end{align*}
with the usual modifications for indices \( p_{i} \) which equal \( \infty \).  

If the measures \( \mu_{i} \) for all \( i = 1, 2, \dots, d \) are clear from context we simply write \( L^{p_{1}, p_{2}, \dots, p_{d}}_{w} \left( X_{1}, X_{2}, \dots, X_{d}\right) \).  If \( X_{i} = \mathbb{R} \) and \(  \mu_{i} \) is Lebesgue measure on \( \mathbb{R} \) for all \( i = 1, 2, \dots, d \), then we simply write \( L^{p_{1}, p_{2}, \dots, p_{d}}_{w} \).  If each \( X_{i} \) is countable and \( \mu_{i} \) is counting measure on \( X_{i} \) we simply write \( \ell^{p_{1}, p_{2}, \dots, p_{d}}_{w} \left( X_{1}, X_{2}, \dots, X_{d}\right) \).

\end{definition}

The mixed norm spaces \( L^{p_{1}, p_{2}, \dots, p_{d}}_{w} \left( X_{1}, X_{2}, \dots, X_{d}, \mu_{1}, \mu_{2}, \dots,  \mu_{d}\right) \) are generalizations of the classical spaces \( L^{p} \) and \(\ell^{p} \), and the proof that \( L^{p} \) and \(\ell^{p} \) are Banach spaces can be extended to the mixed norm spaces (see \cite{oldie}).

The following technical lemma will be useful in later sections.

\begin{lemma} \label{lemma:embedding} If \( p > \frac{2d}{d+s} \) and \( s \geq 0 \) then \( \ell^{2,2}_{v_{s}} \left( \mathbb{Z}^{2d}, \mathbb{Z}^{2d}  \right) \subsetneq \ell^{2,p} \left( \mathbb{Z}^{2d}, \mathbb{Z}^{2d} \right) \).
\end{lemma}

\subsection{Schatten class operators}

\begin{definition}  Suppose \( H \) is a Hilbert space and \( A: H \to H \) is a linear operator.  We say \( A \) is \emph{Schatten \( p \)-class} and write \( A \in \mathcal{I}_{p}(H) \) if \[ \left\| A \right\|_{\mathcal{I}_{p}} = \sup \biggl( \sum_{n \in \mathbb{N}} \left| \langle A f_{n}, g_{n} \rangle \right|^{p} \biggr)^{\frac{1}{p}} < \infty,\] where the supremum is taken over all pairs of orthonormal sequences \( \left\{ f_{n} \right\}_{n \in \mathbb{N}} \),  \( \left\{ g_{n} \right\}_{n \in \mathbb{N}} \) in \( H \).
\end{definition}

Equivalently, an operator is Schatten \( p\)-class if its singular values constitute an \( \ell^{p} \) sequence.  Consequently, trace-class operators are exactly the Schatten 1-class operators and Hilbert-Schmidt operators are the Schatten 2-class operators.  Schatten \( \infty\)-class operators are bounded operators.

\subsection{Gabor Transform}

Suppose \( f:\mathbb{R}^{d} \to \mathbb{C} \) is measurable.  For \( x, \xi \in \mathbb{R}^{d} \) define the translation operator \( T_{x} \) and modulation operator \( M_{\xi} \) by
\[ T_{x}f(t) = f( t -x ) \hspace{2pc} \textrm{ and } \hspace{2pc}  M_{\xi}f(t) = e^{2 \pi i t \cdot \xi} f(t). \]

\begin{definition}  Fix \( \phi \in S(\mathbb{R}^{d}) \).  Given \( f \in S'(\mathbb{R}^{d}) \), the \emph{Gabor transform} of \( f \) with respect to \( \phi \) is 
\[ V_{\phi} f( x, \xi ) = \langle f, M_{\xi} T_{ x} \phi \rangle,  \hspace{2pc} x, \xi \in \mathbb{R}^{d}. \]  The function \( \phi \) is called the \emph{window function} of the Gabor transform.
\end{definition}

The value of \( V_{\phi} f( x, \xi ) \) gives information about the time-frequency content of \( f \) around \( x \) in time and \( \xi \) in frequency.  See \cite{groch} for background and information about the Gabor transform.

\subsubsection{Gabor Frames}

\begin{definition} A \emph{frame} for a Hilbert space \( H \) is a sequence of elements \( \left\{ \phi_{x} \right\}_{x \in \Lambda} \) in \( H \) such that there are \( A, B > 0 \) with 
\[ A \left\| f \right\|^{2} \leq \sum_{x \in \Lambda} \left| \langle f, \phi_{x} \rangle \right|^{2} \leq B \left\| f \right\|^{2}\] for all \( f \in H \).  In this case \( A, B \) are \emph{frame bounds}.  If we can take \( A = B \) then 
\( \left\{ \phi_{x} \right\}_{x \in \Lambda} \)  is a \emph{tight frame}.  A tight frame is \emph{Parseval} if \( A = B = 1 \).  A \emph{Gabor frame} for \( L^{2}(\mathbb{R}^{d}) \) is a sequence \( \left\{ M_{\xi} T_{x} \phi \right\}_{ \left( x, \xi \right) \in \Lambda} \) that is a frame for \( L^{2}( \mathbb{R}^{d}) \). 
\end{definition}

Frames give nonorthogonal expansions of elements of \( H \) in terms of the frame elements.  In particular, if \( \left\{ \phi_{x} \right\}_{x \in \Lambda} \) is a tight frame for \( H \) with frame bound \( B \), we have
\[ f = B^{-1} \sum_{x \in \Lambda} \langle f , \phi_{x} \rangle \phi_{x} \hspace{2pc} \forall f \in H.\] See \cite{lsu} for general background on frames and \cite{groch} for examples and properties of Gabor frames.  In particular, there are tight Gabor frame for \( L^{2}(\mathbb{R}^{d}) \) whose generator \( \phi \) is a nice function, e.g., \( \phi \in C^{\infty}_{c}(\mathbb{R}^{d}) \).  However, the different statements of the Balian-Low Theorem show that the elements of a Gabor frame which offers unique expansions (i.e. a Gabor Riesz basis) necessarily have poor time-frequency localization.

\subsubsection{Wilson Bases}

Wilson bases are orthonormal bases similar to Gabor Riesz bases in that they allow for unique, discrete expansions of the elements of \( L^{2}(\mathbb{R}^{d}) \) in terms of time-frequency ``molecules.''  However, in contrast with Gabor Riesz bases, the elements of a Wilson bases may be well-localized in time and frequency.

For each \(k \in \mathbb{Z}^{d}, n \in \left(  \mathbb{Z}^{+} \right)^{d} \) let \[ \Psi_{k,n}(t) = \psi_{k_{1},n_{1}}(t_{1}) \psi_{k_{2},n_{2}}(t_{2}) \cdots \psi_{k_{d},n_{d}}(t_{d}),\]  where 
\begin{displaymath}
 \psi_{k_{i},n_{i}}(t_{i}) =  \left\{ \begin{array}{ll}
 T_{k_{i}} \psi(t_{i}), & \textrm{ if }  n_{i} = 0,  \\
 \frac{1}{\sqrt{2}} T_{\frac{k_{i}}{2}} \left( M_{n_{i}} + \left( -1 \right)^{k_{i}+n_{i}} M_{-n_{i}} \right) \psi(t_{i}), & \textrm{ if } n_{i} > 0.
 \end{array} \right.
 \end{displaymath}
 
For suitable \( \psi \in L^{2}(\mathbb{R}) \), the sequence \( \left\{ \Psi_{k,n} \right\}_{k \in \mathbb{Z}^{d}, n \in \left(  \mathbb{Z}^{+} \right)^{d}} \) constitutes an orthonormal basis for \( L^{2}(\mathbb{R}^{d}) \).  In this case we call \( \left\{ \Psi_{k,n} \right\}_{k \in \mathbb{Z}^{d}, n \in \left(  \mathbb{Z}^{+} \right)^{d}} \) the \emph{Wilson basis} generated by \( \psi \) (see \cite{groch} for details).
\subsubsection{Modulation Spaces}

Fix \( \phi \in S(\mathbb{R}^{d}) \) and \( p, q \in [1, \infty] \).  Define
%\[ \left\| f \right\|_{M_{w}^{p,q}( \mathbb{R}^{d})} = \left\| V_{g} f \right\|_{L_{w}^{\overbrace{p,p,\dots,p}^{d \textrm{ times}}, \overbrace{q,q,\dots,q}^{d \textrm{ times}}}} \] 
\[ \left\| f \right\|_{M_{w}^{p,q}( \mathbb{R}^{d} )} = \left\| V_{\phi} f \right\|_{L_{w}^{p_{1}, p_{2}, \dots, p_{2d}}}, \] 
where \( p = p_{1} = p_{2} = \dots = p_{d} \) and \( q = p_{d+1} = p_{d+2} = \dots  = p_{2d} \).
Let
\[ M_{w}^{p,q}(\mathbb{R}^{d}) = \left\{ f \in S'(\mathbb{R}^{d}) : \left\| f \right\|_{M_{w}^{p,q}( \mathbb{R}^{d})} < \infty \right\}. \]  Each \( M_{w}^{p,q}(\mathbb{R}^{d}) \) is a \emph{modulation space}. For \( w = 1 \) we write \( M_{w}^{p,q}( \mathbb{R}^{d}) = M^{p,q}( \mathbb{R}^{d}) \).

The modulation space \( M_{w}^{p,q}\left(\mathbb{R}^{d}\right) \) consists of functions with a particular time-frequency decay controlled by the parameters \( p,q \) and weight \( w \).

\subsection{Integral operators and Pseudodifferential Operators}

An operator \( A \) of the form 
\[ Af(t) = \int_{\mathbb{R}^{d}} k(t,y) f(y) \, \ud y\] is an \emph{integral operator}.  The function \( k \) is the \emph{kernel} of \( A \).

A \emph{pseudodifferential operator} with \emph{Kohn-Nirenberg symbol} \( \tau \) is an operator having the form
\[ K_{\tau} f(t) = \iint_{\mathbb{R}^{2d}} \hat{\tau} \left( \xi, x \right)  M_{\xi} T_{-x} f(t) \, \ud x \, \ud \xi. \] 

Suitable \(  K_{\tau} \) can be realized as integral operators.  In particular, if we let \( \mathcal{F}_{2} \) denote the partial Fourier transform on the last \( d \) variables of a function of \( 2d \) variables, i.e.
\[ \left( \mathcal{F}_{2} F \right)(x,w) = \int_{\mathbb{R}^{d}}  F(x,y)\, e^{-2 \pi i y \cdot w} \, \ud y  \hspace{2pc} \textrm{ for all } x, w \in \mathbb{R}^{d} ,\] then \( K_{\tau}  \) is an integral operator with kernel \( k = \mathcal{F}_{2}^{-1} \tau \circ N \), where \( N(x,y) = \left( x, x-y \right) \) for \( x, y \in \mathbb{R}^{d} \).

%\subsubsection{Schatten \(p\)-class Pseudodifferential Operators}

In general, the time-frequency properties of the symbol of a pseudodifferential operator determine if the operator is Schatten \(p\)-class.  Results of this type can be found in  \cite{modbounded}, \cite{heilpdo}, \cite{roch1}, \cite{sjo}, and \cite{toftpdo}.  The strongest of these results is found in \cite{modbounded}, in which the authors obtain estimates on the singular values of pseudodifferential operators.  The following theorem is a special case of Theorem 1.2 in \cite{modbounded}.

\begin{theorem} \label{thm:heilgrochpclass} Suppose \( A \) is a pseudodifferential operator with Kohn-\\Nirenberg symbol \( \tau \).  If \( \tau \in M^{2,2}_{v_{s}}(\mathbb{R}^{2d}) \) with \( p > \frac{2d}{d+s} \) and \( s \geq 0 \), then \( A \in \mathcal{I}_{p}\left(L^{2 }(\mathbb{R}^{d}) \right) \).
\end{theorem}

\section{A Schatten Class result for integral operators}

In this section we find a general condition on the kernel of an integral operator which ensures the operator is Schatten \(p \)-class.

\begin{lemma} \label{lemma:pclassG}  Assume \( \left\{ f_{j} \right\}_{j \in \mathbb{N}}, \left\{ g_{j} \right\}_{j \in \mathbb{N}} \) are orthonormal sequences in \( L^{2}(\mathbb{R}^{d}) \).  Suppose \( \left\{ \phi_{n} \right\}_{n \in \Lambda} \) is a Parseval frame for \( L^{2}(\mathbb{R}^{d}) \). For \( G \in L^{2,p}(\mathbb{R}^{d}, \Lambda) \) define 
\[ T(G) = \left\{ \sum_{n \in \Lambda} \langle f_{j}, \phi_{n}  \rangle \langle G(\cdot, n), g_{j}\rangle \right\}_{j \in \mathbb{N}}. \]  Then for all \( p \in [1,2] \), \( T: L^{2,p}(\mathbb{R}^{d}, \Lambda) \to \ell^{p}\left( \mathbb{N} \right) \) is bounded  with \( \left\| T \right\| \leq 1 \).
\end{lemma}
\begin{proof}  Since \( \left\{ \phi_{n} \right\}_{n \in \Lambda} \) has frame bounds \( A = B = 1\), we have \( \left\| \phi_{n} \right\|_{L^{2}(\mathbb{R}^{d})} \leq 1 \) for all \( n \in \Lambda\).  Therefore
\begin{align*}
\left\| T(G) \right\|_{\ell^{1}} & =  \sum_{j \in \mathbb{N}} \, \biggl|\sum_{n \in \Lambda} \langle f_{j}, \phi_{n}  \rangle \langle G(\cdot, n), g_{j}\rangle \biggr| \\[2 \jot]
& \leq  \sum_{n \in \Lambda} \sum_{j\in \mathbb{N}} \left| \langle f_{j}, \phi_{n} \rangle  \right| \left| \langle G(\cdot, n), g_{j}\rangle \right|  \\
& \leq  \sum_{n \in \Lambda} \biggl( \sum_{j\in \mathbb{N}} \left| \langle f_{j}, \phi_{n} \rangle  \right|^{2} \biggr)^{\frac{1}{2}} \biggl( \sum_{j\in \mathbb{N}} \left| \langle G(\cdot, n) , g_{j} \rangle  \right|^{2} \biggr)^{\frac{1}{2}} \\[2 \jot]
& \leq   \sum_{n \in \Lambda}   \left\| \phi_{n} \right\|_{L^{2}(\mathbb{R}^{d})} \left\| G(\cdot, n) \right\|_{L^{2}(\mathbb{R}^{d})} \\
& \leq    \left\| G \right\|_{L^{2,1}(\mathbb{R}^{d}, \Lambda)}
\end{align*}
and
\begin{align*}
\left\| T(G) \right\|_{\ell^{2}} & = \biggl( \sum_{j \in \mathbb{N}} \, \biggl|  \sum_{n \in \Lambda} \langle f_{j}, \phi_{n} \rangle \langle G(\cdot, n), g_{j}\rangle  \biggr|^{2} \biggr)^{\frac{1}{2}} \\
& \leq \biggl( \sum_{j \in \mathbb{N}}   \biggl(  \sum_{n \in \Lambda} \left| \langle f_{j}, \phi_{n} \rangle \right|^{2}   \biggr)  \biggl(  \sum_{n \in \Lambda} \left| \langle G(\cdot, n), g_{j}\rangle \right|^{2}   \biggr)   \biggr)^{\frac{1}{2}} \\
& =  \biggl(     \sum_{n \in \Lambda}  \sum_{j \in \mathbb{N}} \left| \langle g_{j}, G(\cdot, n) \rangle \right|^{2}   \biggr)^{\frac{1}{2}} \\
& \leq  \biggl(   \sum_{n \in \Lambda}   \left\|  G(\cdot, n)  \right\|_{L^{2}(\mathbb{R}^{d})}^{2}     \biggr)^{\frac{1}{2}} \\[1 \jot]
& =   \left\| G \right\|_{L^{2,2}(\mathbb{R}^{d}, \Lambda)}.
\end{align*}
Hence the theorem holds for \( p = 1 \) and \( p = 2 \).  The Riesz-Thorin Interpolation Theorem gives the result for \( p \in (1,2) \). 
\end{proof}

\begin{theorem} \label{thm:modsum} Suppose \( \left\{ \phi_{m} \right\}_{m \in \Lambda} \) is a Parseval frame for \( L^{2}(\mathbb{R}^{d}) \).  Define \( \Phi_{m,n}(t,y) = \phi_{m}(t) \overline{\phi_{n}(y)}  \).  If \( A \) is an integral operator with kernel \( k \) then for all \( p \in [1,2] \)
\[ \left\| A \right\|_{\mathcal{I}_{p}} \leq \biggl( \sum_{n \in \Lambda} \biggl( \sum_{m \in \Lambda} \left| \langle k, \Phi_{m,n} \rangle \right|^{2} \biggr)^{\frac{p}{2}}  \biggr)^{\frac{1}{p} }. \] 
\end{theorem}
\begin{proof} Suppose \( \left\{ f_{j} \right\}_{j \in \mathbb{N}}, \left\{ g_{j} \right\}_{j \in \mathbb{N}} \) are orthonormal sequences in \( L^{2}(\mathbb{R}^{d}) \).  Let \( G(y,n) = A\phi_{n}(y) \).  Notice that \( \langle \phi_{n}, A^{*} g_{j} \rangle = \langle G(\cdot,n), g_{j} \rangle \).  Expanding \( f_{j} \) with the frame \( \left\{ \phi_{m} \right\}_{m \in \Lambda} \) and using the previous lemma, we have 
\begin{align*}
\biggl( \sum_{j \in \mathbb{N}} \left| \langle A f_{j}, g_{j} \rangle \right|^{p} \biggr)^{\frac{1}{p}} 
&= \biggl( \sum_{j \in \mathbb{N}} \left| \langle  f_{j}, A^{*} g_{j} \rangle \right|^{p} \biggr)^{\frac{1}{p}} \\
& =  \biggl( \sum_{j \in \mathbb{N}} \biggl|  \sum_{n \in \Lambda} \langle  f_{j}, \phi_{n} \rangle \langle G(\cdot, n), g_{j} \rangle  \biggr|^{p} \biggr)^{\frac{1}{p}} \\[1 \jot]
& \leq \left\| G \right\|_{L^{2,p}(\mathbb{R}^{d}, \Lambda)} \\
& =  \biggl( \sum_{n \in \Lambda} \left\| A \phi_{n} \right\|^{p}_{L^{2}(\mathbb{R}^{d})}  \biggr)^{\frac{1}{p} }\\
& = \biggl( \sum_{n \in \Lambda} \biggl( \sum_{m \in \Lambda} \left| \langle A \phi_{n}, \phi_{m} \rangle \right|^{2} \biggr)^{\frac{p}{2}}  \biggr)^{\frac{1}{p} }\\
& = \biggl( \sum_{n \in \Lambda} \biggl( \sum_{m \in \Lambda} \left| \langle k, \Phi_{m,n} \rangle \right|^{2} \biggr)^{\frac{p}{2}}  \biggr)^{\frac{1}{p} }.
\end{align*}
Taking the supremum of \( \left( \sum_{j \in \mathbb{N}} \left| \langle A f_{j}, g_{j} \rangle \right|^{p} \right)^{\frac{1}{p}}  \) over all such orthonormal sequences \( \left\{ f_{j} \right\}_{j \in \mathbb{N}}, \left\{ g_{j} \right\}_{j \in \mathbb{N}} \) gives the result. \qedhere

\end{proof}

The proofs of Lemma \ref{lemma:pclassG} and Theorem \ref{thm:modsum} can be generalized to prove Theorem \ref{thm:pclassframe}.

\section{Mixed Modulation Spaces}

In this section we introduce a generalization of the modulation spaces \( M_{w}^{p,q}(\mathbb{R}^{d}) \). Throughout this section, we assume \( c \) is a permutation of the set \( \left\{ 1,2,\dots, 2d\right\}  \).  To simplify some notation we identify  \(c \) with the bijection \( \mathfrak{c}: \mathbb{R}^{2d} \to \mathbb{R}^{2d} \) given by \( \mathfrak{c}(x_{1},\dots x_{2d}) = (x_{c(1)}, \dots, x_{c(2d)}) \).

\begin{definition}  Suppose \( \phi \in S(\mathbb{R}^{d}) \) and \( c \) is a permutation of \( \left\{ 1,2,\dots, 2d\right\}  \) corresponding to the map \( \mathfrak{c} \).  Let \( M(c)_{w}^{p_{1}, p_{2}, \dots, p_{2d}} \) be the \emph{mixed modulation space} consisting of all \(  f \in S'(\mathbb{R}^{d}) \) for which
\[ \left\| f \right\|_{M(c)_{w}^{p_{1}, p_{2}, \dots, p_{2d}}} = \left\| V_{\phi} f \circ \mathfrak{c} \right\|_{L_{w}^{p_{1}, p_{2}, \dots, p_{2d}}} < \infty. \]  When \( w = 1 \) we write \( M(c)_{w}^{p_{1}, p_{2}, \dots, p_{2d}} = M(c)^{p_{1}, p_{2}, \dots, p_{2d}} \).
\end{definition}

Notice that if \( c \) is the identity permutation  and \( p = p_{1} = p_{2} = \dots = p_{d} \) and \( q = p_{d+1} = \dots = p_{2d} \) then \( M(c)_{w}^{p_{1}, p_{2}, \dots, p_{2d}} = M^{p,q}_{w}(\mathbb{R}^{d})  \). Hence the mixed modulation spaces are indeed generalizations of modulation spaces. Also notice that if  \( p = p_{1} = p_{2} = \dots = p_{d} = p_{d+1} = \dots = p_{2d} \) then \( M(c)_{v_{s}}^{p_{1}, p_{2}, \dots, p_{2d}} = M^{p,p}_{v_{s}}(\mathbb{R}^{d})  \) for any permutation \( c \).

The most interesting properties of modulation spaces carry over to the mixed modulation spaces.  As the proofs are basic generalizations of the proofs for modulation spaces, we state these properties without proof.  See \cite{groch} for a detailed account of the properties of modulation spaces.

\begin{definition} Suppose \( c \) is a permutation of \( \left\{ 1,2,\dots, 2d\right\}  \).  For each \( x  \in \mathbb{R}^{2d} \) let \( \pi_{x} = M_{\left( x_{d+1}, \dots, x_{2d} \right)} T_{\left( x_{1}, \dots, x_{d} \right)} \).  For measurable \( \psi:\mathbb{R}^{d} \to \mathbb{C} \) define an operator \( \Upsilon_{\psi} \) by
\[ \Upsilon_{\psi} F(t) = \int_{\mathbb{R}^{2d}} F(x) \, \pi_{\mathfrak{c}(x)} \psi(t) \, \ud x. \]
\end{definition}

\begin{theorem} \label{thm:prop} Suppose \( \psi, \gamma \in M(c)^{1,\dots,1}_{v} \).  
\begin{itemize}
\item[(a)]  For any \( f \in M(c)_{w}^{p_{1}, p_{2}, \dots, p_{2d}} \), we have \( \Upsilon_{\psi} \left( V_{\gamma} f                                                        \circ \mathfrak{c} \right) = \langle \psi, \gamma \rangle f \). \smallskip
\item[(b)]  \( ||| f ||| = \left\| V_{\psi} f \circ \mathfrak{c} \right\|_{L_{w}^{p_{1}, p_{2}, \dots, p_{2d}}} \) is an equivalent norm on \( M(c)_{w}^{p_{1}, p_{2}, \dots, p_{2d}} \).
\end{itemize}
 
\end{theorem}  

Theorem \ref{thm:prop}(b) implies that the definition of the mixed modulation spaces is independent of the choice of window \( \phi \in S(\mathbb{R}^{d}) \), with different windows \( \phi \) giving equivalent norms.  Furthermore, this fact also holds for \( \phi \) in the larger space \( M(c)^{1,\dots,1}_{v}  \).  Theorem \ref{thm:prop}(a) shows that for Gabor window functions in \( M(c)^{1,\dots,1}_{v}  \), there is an inversion formula valid on each \( M(c)_{w}^{p_{1}, p_{2}, \dots, p_{2d}} \).

\begin{corollary} For any \( p_{1}, p_{2}, \dots, p_{2d} \in [1, \infty] \), \(  M(c)_{w}^{p_{1}, p_{2}, \dots, p_{2d}} \) is a Banach space.
\end{corollary}

\begin{theorem}  If \( p_{1}, p_{2}, \dots, p_{2d} \in [1, \infty) \) then \(  M(c)_{\frac{1}{w}}^{p'_{1}, p'_{2}, \dots, p'_{2d}} \) is the dual space of \(  M(c)_{w}^{p_{1}, p_{2}, \dots, p_{2d}} \), where \( p'_{i} \in [1, \infty] \) satisfies \( \frac{1}{p_{i}} + \frac{1}{p'_{i}} = 1 \).
\end{theorem}

The next theorem states that if the window function is nice then a Gabor frame for \( L^{2}(\mathbb{R}^{d}) \) gives bounded decompositions for all mixed modulation spaces.

\begin{theorem}  \label{thm:banachframe} Suppose \(  p_{1}, p_{2}, \dots, p_{2d} \in [1, \infty] \) and \( \psi \in M(c)^{1,\dots,1}_{v} \).  Further suppose that \( \left\{ \pi_{\alpha n} \psi \right\}_{n \in \mathbb{Z}^{2d}} \) is a frame for \( L^{2}(\mathbb{R}^{d}) \) with dual frame \( \left\{ \pi_{\alpha n} \gamma \right\}_{n \in \mathbb{Z}^{2d}} \). Then
\begin{itemize}
\item[(a)] \( \left\{ \pi_{\alpha n} \psi \right\}_{n \in \mathbb{Z}^{2d}} \) is a Banach frame for \(  M(c)_{w}^{p_{1}, p_{2}, \dots, p_{2d}} \) and there exist \(0 < A \leq B < \infty \) independent of \(   p_{1}, p_{2}, \dots, p_{2d} \) with \[ A \left\| f \right\|_{ M(c)_{w}^{p_{1}, p_{2}, \dots, p_{2d}}} \leq \left\| V_{\psi} f \circ \mathfrak{c} \big|_{\alpha \mathbb{Z}^{2d}} \right\|_{\ell_{w}^{p_{1}, p_{2}, \dots, p_{2d}}} \leq B \left\| f \right\|_{ M(c)_{w}^{p_{1}, p_{2}, \dots, p_{2d}}}, \] for all \( f \in M(c)_{w}^{p_{1}, p_{2}, \dots, p_{2d}} \).
\item[(b)]  If \( p_{1}, p_{2}, \dots, p_{2d} \in [1, \infty) \) then \begin{displaymath} f = \sum_{m \in \mathbb{Z}^{2d}} \langle f, \pi_{\alpha m} \psi \rangle \, \pi_{\alpha m} \gamma =  \sum_{m \in \mathbb{Z}^{2d}} \langle f, \pi_{\alpha m} \gamma \rangle \, \pi_{\alpha m} \psi \end{displaymath} for all \( f \in M(c)_{w}^{p_{1}, p_{2}, \dots, p_{2d}} \) with unconditional convergence in \( M(c)_{w}^{p_{1}, \dots, p_{2d}} \). \smallskip
\item[(c)]  If \( p_{1}, p_{2}, \dots, p_{2d} \in [1, \infty] \) then  \begin{displaymath} f = \sum_{m \in \mathbb{Z}^{2d}} \langle f, \pi_{\alpha m} \psi \rangle \pi_{\alpha m} \gamma =  \sum_{m \in \mathbb{Z}^{2d}} \langle f, \pi_{\alpha m} \gamma \rangle \pi_{\alpha m} \psi \end{displaymath} for all \( f \in M(c)_{w}^{p_{1}, p_{2}, \dots, p_{2d}} \) with weak* convergence in \( M(c)^{\infty, \dots, \infty}_{\frac{1}{v}} \).
\end{itemize}
\end{theorem}

Theorem \ref{thm:banachframe} can be used to prove embeddings among the mixed modulation spaces.

\begin{lemma} \label{lemma:modembed} If \( s \geq t \) and \( p_{i}, r_{i} \in [1, \infty] \) with \( p_{i} \leq r_{i} \) for all \( i = 1,2, \dots, 2d \) then \(M(c)_{v_{s}}^{p_{1},p_{2},\dots , p_{2d}} \subset  M(c)_{v_{t}}^{r_{1},r_{2},\dots ,r_{2d}}  \). 
\end{lemma}

The following theorem states that Wilson bases are bases for the mixed modulation spaces.

 \begin{theorem} \label{thm:wilsonbasis}  Let \( v:\mathbb{R}^{2d} \to [0, \infty) \) be a weight and \( w \) a \( v \)-moderate weight.  Define \(\tilde{v}(t) = \max \left\{ v(t,0, \dots, 0), v(0,t,0, \dots, 0) , \dots, v(0, \dots, 0,t) \right\} \) for each \( t \in \mathbb{R} \) .  Assume \( \psi  \in M^{1,1}_{\tilde{v}\otimes \tilde{v}}(\mathbb{R})\) generates an orthonormal Wilson basis \( \left\{ \Psi_{k,n} \right\}_{n \in \left( \mathbb{Z}^{+} \right)^{d} , k \in \mathbb{Z}^{d}} \) for \( L^{2}(\mathbb{R}^{d}) \).  Then  \( \left\{ \Psi_{k,n} \right\}_{n \in \left( \mathbb{Z}^{+} \right)^{d} , k \in \mathbb{Z}^{d}} \) is an unconditional basis for
 \(M(c)_{w}^{p_{1},p_{2}, \dots, p_{2d}}  \) for each \( p_{1},p_{2}, \dots, p_{2d} \in [1, \infty) \).
 \end{theorem}

 \begin{corollary} \label{cor:isothm} Let \( X_{1} = X_{2} = \dots = X_{d} = \mathbb{Z} \) and \( X_{d+1} = X_{d+2} = \dots = X_{2d} = \mathbb{Z}^{+} \).  Then the map 
 \[ f \to \left\{ \langle f, \Psi_{(n_{c(1)},n_{ c(2)}, \dots, n_{c(d)}),(n_{c(d+1)}, \dots, n_{c(2d)})} \rangle \right\}_{n_{1} \in X_{c^{-1}(1)}, n_{2} \in X_{c^{-1}(2)}, \dots, n_{2d} \in X_{c^{-1}(2d)}}\] is an isomorphism from \( M(c)_{w}^{p_{1},p_{2}, \dots, p_{2d}} \) to \( \ell_{w}^{p_{1},p_{2}, \dots, p_{2d}} \left( X_{c^{-1}(1)}, \dots , X_{c^{-1}(2d)} \right) \).
 \end{corollary}

\section{Pseudodifferential Operators and Schatten classes}
In this section we will use Theorem \ref{thm:modsum} to find conditions on the kernel and Kohn-Nirenberg symbol of a pseudodifferential operator that guarantee the operator is Schatten \( p \)-class. 

We are particularly interested in permutations \( c \) of \( \left\{ 1,2, \dots, 4d \right\}  \) satisfying  the following:
\begin{itemize}
\item[(a)] \( c \) maps \( \left\{ 1,2,\dots, d, 2d+1, 2d+2, \dots, 3d \right\}\) to  \( \left\{ 1,2,\dots, 2d \right\} \)  bijectively and
\item[(b)] \( c \) maps  \( \left\{d+1,d+2, \dots, 2d, 3d+1, 3d+2, \dots, 4d\right\} \) to \( \left\{ 2d+1, \dots, 4d \right\} \)  bijectively.
\end{itemize}  We call such permutations \emph{slice permutations} because they relate nicely to the slice analysis of Section 3.

\begin{corollary} \label{cor:kernelmod}    Assume \( c \) is a slice permutation.  Let \( 2 = p_{1} = p_{2} = \dots = p_{2d} \) and \( p = p_{2d+1} = \dots = p_{4d} \). If \( p \in [1,2] \), \( k \in M(c)^{p_{1},\dots, p_{4d}} \) and \( A \) is an integral operator with kernel \( k \), then \( A \in \mathcal{I}_{p}(L^{2}(\mathbb{R}^{d})) \).
\end{corollary}
\begin{proof} Let \( \left\{ \pi_{\alpha m} \phi \right\}_{m \in \mathbb{Z}^{2d}} = \left\{ \phi_{m} \right\}_{m \in \mathbb{Z}^{2d}} \) be a Parseval Gabor frame for \( L^{2}(\mathbb{R}^{d}) \) with \( \phi \in M^{1,1}(\mathbb{R}^{d}) \) and let \( \Phi(t,y) = \phi(t) \overline{\phi(y)} \).  Then \( \Phi \in M^{1,1}(\mathbb{R}^{2d}) \).  Let \( \Phi_{m, n}(t,y) = \phi_{m}(t) \overline{\phi_{n}(y)} \).  By Lemma 3.2 in \cite{singval}, \( \left\{  \Phi_{m, n} \right\}_{m, n \in \mathbb{Z}^{2d}} \) is a Parseval frame for \( L^{2}(\mathbb{R}^{2d}) \).  
For \( m_{1}, m_{2}, n_{1},n_{2} \in \mathbb{Z}^{d} \), with \( m = (m_{1}, m_{2} ) \) and  \( n = (n_{1}, n_{2} ) \),  we have 
\[ \langle k ,  \Phi_{m,n} \rangle = V_{\Phi}k( \alpha m_{1}, \alpha n_{1}, \alpha m_{2}, \alpha n_{2}). \] For each slice permutation \(c \), we see that
\begin{align*}
& \biggl( \sum_{n \in \mathbb{Z}^{2d}} \biggl( \sum_{m \in \mathbb{Z}^{2d}} \left| \langle k, \Phi_{m,n} \rangle \right|^{2} \biggr)^{\frac{p}{2}}\biggr)^{\frac{1}{p}} \\
& \indent \indent \indent \indent = \biggl( \sum_{n_{1}, n_{2} \in \mathbb{Z}^{d}} \biggl( \sum_{m_{1}, m_{2} \in \mathbb{Z}^{d}} \left| V_{\Phi}k( \alpha m_{1}, \alpha n_{1}, \alpha m_{2}, \alpha n_{2}) \right|^{2} \biggr)^{\frac{p}{2}}\biggr)^{\frac{1}{p}} \\
& \indent \indent \indent \indent = \biggl( \sum_{n_{1}, n_{2} \in \mathbb{Z}^{d}} \biggl( \sum_{m_{1}, m_{2} \in \mathbb{Z}^{d}} \left| V_{\Phi}k( \mathfrak{c}(\alpha m_{1}, \alpha m_{2}, \alpha n_{1}, \alpha n_{2})) \right|^{2} \biggr)^{\frac{p}{2}}\biggr)^{\frac{1}{p}} \\[2 \jot]
& \indent \indent \indent \indent \leq B  \left\| k \right\|_{M(c)^{p_{1},p_{2}, \dots, p_{4d}}},
\end{align*}
where \( B \) is the constant ensured by Theorem \ref{thm:banachframe}(a).  Hence if \( k \in M(c)^{p_{1}, \dots, p_{4d}}  \) then, by Theorem \ref{thm:modsum}, \( A \in \mathcal{I}_{p}(L^{2}(\mathbb{R}^{d})) \).
\end{proof}

We can extend Corollary \ref{cor:kernelmod} to conditions on the symbol of a pseudodifferential operator.  

%\begin{corollary} \label{cor:forcomparison} Assume \( c: \left\{ 1,2, \dots, 4d \right\} \to  \left\{ %1,2, \dots, 4d \right\} \) is a slice permutation.  Let \( A \) be a pseudodifferential operator with %kernel \( k \), Kohn-Nirenberg symbol \( \tau \) and Weyl symbol \( \sigma \).  Assume  \( p \in [1,2] %\) and set \( 2 = p_{1} = \dots = p_{2d} \) and \( p = p_{2d+1} = \dots = p_{4d} \). If one of \( k, %\tau, \sigma  \) lies in \( M(c)^{p_{1},p_{2},\dots, p_{4d}} \) then do are the others.  In this case %\( A \in \mathcal{I}_{p}\left( L^{2}(\mathbb{R}^{d}) \right) \).
%\end{corollary}
\begin{theorem} \label{thm:maintfspec} Let \( A \) be a pseudodifferential operator with kernel \( k \) and Kohn-Nirenberg symbol \( \tau \).  Assume  \( p \in [1,2] \) and set \( 2 = p_{1} = \dots = p_{2d} \) and \( p = p_{2d+1} = \dots = p_{4d} \). If \( c \) is a slice permutation and one of \( k, \tau  \) lies in \( M(c)^{p_{1},p_{2},\dots, p_{4d}} \) then so does the other.  In this case \( A  \in \mathcal{I}_{p}(L^{2}(\mathbb{R}^{d})) \).
\end{theorem}
\begin{proof} Using the fact that \( k = \mathcal{F}_{2}^{-1} \tau \circ N \) we can show that \( \left|\langle k , M_{(z,t)} T_{(x,y)} \Phi \rangle \right|  = \left|\langle \tau ,M_{(z+t,-y)} T_{(x,-t)} \mathcal{F}_{2}( \Phi \circ N^{-1}) \right| \) for all \( x,y,z,t \in \mathbb{R}^{d} \).  Hence
\begin{align*}
&  \left\| k \right\|_{M(c)^{p_{1},p_{2},\dots, p_{4d}}} \\
& \indent =  \left( \iint \left( \iint  \left|\langle k , M_{(z,t)} T_{(x,y)} \Phi \rangle \right|^{2} \, \ud x \, \ud z \right)^{\frac{p}{2}} \, \ud y \, \ud t\right)^{\frac{1}{p}} \\
& \indent = \left( \iint \left( \iint  \left|\langle \tau ,M_{(z+t,-y)} T_{(x,-t)} \mathcal{F}_{2}( \Phi \circ N^{-1}) \right|^{2} \, \ud x \, \ud z \right)^{\frac{p}{2}} \, \ud y \, \ud t \right)^{\frac{1}{p}} \\
& \indent = \left( \iint \left( \iint  \left|\langle \tau , M_{(z,y)} T_{(x,t)} \mathcal{F}_{2}( \Phi \circ N^{-1}) \rangle \right|^{2} \, \ud x \, \ud z \right)^{\frac{p}{2}} \, \ud y \, \ud t \right)^{\frac{1}{p}} \\[2 \jot]
& \indent \asymp \left\| \tau \right\|_{M(c)^{p_{1},p_{2},\dots, p_{4d}}}. \qedhere
\end{align*}

\end{proof}  

%Corollary \ref{cor:forcomparison} 
Theorem \ref{thm:maintfspec} is stronger than the previously known Theorem \ref{thm:heilgrochpclass}, as the following lemma shows.

\begin{lemma} Let \( c \) be a slice permutation and let \( 2 = p_{1} = \dots = p_{2d} \), \( p = p_{2d+1} = \dots = p_{4d} \).  If \( p > \frac{2d}{d+s} \) with \( s \geq 0 \), then  \( M^{2,2}_{v_{s}}(\mathbb{R}^{2d}) \subsetneq M(c)^{p_{1},p_{2},\dots, p_{4d}} \).
\end{lemma}  
\begin{proof}  
Let \( X_{1} = X_{2} = \dots = X_{2d} = \mathbb{Z} \) and \( X_{2d+1} = X_{2d+2} = \dots = X_{4d} = \mathbb{Z}^{+} \) and define  \[ S(f) = \left\{ \langle f, \Psi_{(n_{c(1)}, \dots, n_{c(2d)}), (n_{c(2d+1)}, \dots, n_{c(4d)})} \rangle \right\}_{n_{1} \in X_{c^{-1}(1)}, n_{2} \in X_{c^{-1}(2)}, \dots, n_{4d} \in X_{c^{-1}(4d)}}.\]  Since \( M^{2,2}_{v_{s}}(\mathbb{R}^{2d}) = M(c)_{v_{s}}^{2,2,\dots, 2} \),  Corollary \ref{cor:isothm} implies that \[ S: M^{2,2}_{v_{s}}(\mathbb{R}^{2d}) \to \ell_{v_{s}}^{2,2, \dots, 2} \left( X_{c^{-1}(1)}, \dots, X_{c^{-1}(4d)} \right) \]  and \[S: M(c)^{p_{1},p_{2},\dots, p_{4d}} \to \ell^{p_{1},p_{2}, \dots, p_{4d}} \left( X_{c^{-1}(1)}, \dots, X_{c^{-1}(4d)} \right) \] are isomorphisms.   Furthermore, by Lemma \ref{lemma:embedding}, we have \[ \ell_{v_{s}}^{2,2, \dots, 2} \left( X_{c^{-1}(1)}, \dots, X_{c^{-1}(4d)} \right) \subsetneq \ell^{p_{1},p_{2}, \dots, p_{4d}} \left( X_{c^{-1}(1)}, \dots, X_{c^{-1}(4d)} \right) \] for \( p > \frac{2d}{d+s} \) with \( s \geq 0 \).  Hence we obtain the following diagram \\

\xymatrix{ M^{2,2}_{v_{s}}(\mathbb{R}^{2d}) \ar[d]_{S} &   M(c)^{p_{1},p_{2},\dots, p_{4d}} \ar[d]_{S} \\ 
\ell_{v_{s}}^{2,2, \dots, 2} \left( X_{c^{-1}(1)}, \dots, X_{c^{-1}(4d)} \right) \ar@{^{(}->}[r] & \ell^{p_{1},p_{2}, \dots, p_{4d}} \left( X_{c^{-1}(1)}, \dots, X_{c^{-1}(4d)} \right).}

Since \( S \) is an isomorphism, the result follows.
\end{proof}

By Lemma \ref{lemma:modembed}, increasing any one of the exponent parameters \( p_{1},  \dots, p_{4d}\) or decreasing the weight parameter \( s \) yields a mixed modulation space larger than \( M(c)_{v_{s}}^{p_{1},p_{2},\dots , p_{4d}} \).  The next theorem shows Theorem \ref{thm:maintfspec} is sharp in the following sense: increasing the exponent parameters or decreasing the weight parameter of the mixed modulation space in Theorem \ref{thm:maintfspec} gives a larger mixed modulation space, but pseudodifferential operators with kernels or Kohn-Nirenberg symbols in this larger space need not be Schatten class.

\begin{theorem} Assume \( s \leq 0 \), \( p_{1} ,  \dots  , p_{2d} \in [2, \infty] \), \(  p_{2d+1} , \dots , p_{4d} \in [p, \infty] \) and \( c \) is a slice permutation. Assume at least one of the following is true:
\begin{itemize}
\item[(a)] \(s < 0 \).
\item[(b)]  At least one of \(  p_{1} ,  \dots  , p_{2d} \) is larger than 2.
\item[(c)] At least one of \(  p_{2d+1} , \dots , p_{4d} \) is larger than \( p \).
\end{itemize} 
If \( 1 \leq p \leq 2 \) then there are pseudodifferential operators with kernels in \( M(c)_{v_{s}}^{p_{1},p_{2},\dots, p_{4d}} \) and pseudodifferential operators with Kohn-Nirenberg symbols in \( M(c)_{v_{s}}^{p_{1},p_{2},\dots, p_{4d}} \) that are not in \( \mathcal{I}_{p}(L^{2}(\mathbb{R}^{d})) \).
\end{theorem}
\begin{proof}  Suppose \( A \) is a pseudodifferential operator with kernel \( k \) and Kohn-Nirenberg symbol \( \tau \).  Since \[ \left|\langle k , M_{(z,t)} T_{(x,y)} \Phi \rangle \right|  = \left|\langle \tau ,M_{(z+t,-y)} T_{(x,-t)} \mathcal{F}_{2}( \Phi \circ N^{-1}) \right|, \] it follows that for each slice permutation \( c \), there is a slice permutation \( \tilde{c} \) with \(  \left\| \tau \right\|_{M(c)^{p_{1},p_{2},\dots, p_{4d}}} \asymp \left\| k \right\|_{M(\tilde{c})^{p_{1},p_{2},\dots, p_{4d}}} \).  Hence it suffices to show that for each slice permutation \( c \), there are pseudodifferential operators with kernels in  \( M(c)_{v_{s}}^{p_{1},p_{2},\dots, p_{4d}} \) that are not Schatten \(p \)-class.

To avoid complicated notation, we prove the theorem only for the permutation \( \mathfrak{c}(x_{1}, \dots, x_{4d})= ( x_{d+1}, \dots, x_{2d}, x_{3d+1}, \dots, x_{4d}, x_{1}, \dots, x_{d}, x_{2d+1}, \dots, x_{3d})\).  The result is proven similarly for other slice permutations.

In the case that (a) or (b) holds, we can adapt some of the arguments in \cite{counterex} to complete the proof.  In particular, if \( k(t,y) = k_{1}(t) k_{2}(y) \) is the kernel of an integral operator \( A \), then \( Af = \langle f, \overline{k_{2}} \rangle k_{1} \).  Hence if  \( k_{1} \notin L^{2}(\mathbb{R}^{d}) \), then \( A \) does not map into \( L^{2}(\mathbb{R}^{d}) \), and if  \( k_{2} \notin L^{2}(\mathbb{R}^{d}) \), then \( A:L^{2}(\mathbb{R}^{d}) \to L^{2}(\mathbb{R}^{d}) \) is not bounded.  Let \( c' \) be the permutation with associated bijection \( \mathfrak{c}'(n_{1}, \dots, n_{2d}) = (n_{d+1}, \dots, n_{2d}, n_{1}, \dots n_{d}) \).  If (a) holds, choose \( k_{1} \in M^{2,2}_{v_{s}}(\mathbb{R}^{d}) \setminus L^{2}(\mathbb{R}^{d}) \) and \( k_{2} \in  M^{p,p}(\mathbb{R}^{d}) \).  If (b) holds, choose \( k_{1} \in M(c')^{p_{1}, \dots, p_{2d}} \setminus L^{2}(\mathbb{R}^{d}) \) and \( k_{2} \in  M(c')^{p_{2d+1}, \dots, p_{4d}} \).  In either case \( k(t,y) = k_{1}(t) k_{2}(y) \in M(c)_{v_{s}}^{p_{1},p_{2},\dots, p_{4d}} \), but the integral operator with kernel \( k \) is not a bounded operator on \( L^{2}(\mathbb{R}^{d} ) \).

Hence we assume (c) is true. Choose \( \lambda \in \ell^{p_{2d+1}, \dots, p_{3d},p_{3d+1}, \dots p_{4d}} ((\mathbb{Z}^{+})^{d},\mathbb{Z}^{d}  ) \setminus \ell^{p, p} ((\mathbb{Z}^{+})^{d},\mathbb{Z}^{d}  ) \). Assume \( \left\{ \psi_{j,l} \right\}_{j \in \mathbb{Z}^{d}, l \in (\mathbb{Z}^{+})^{d}} \) is a Wilson basis for \( L^{2}(\mathbb{R}^{d}) \) generated by \( \psi \in M^{1,1}(\mathbb{R}) \).  Then \[ \left\{ \Psi_{(j_{1},j_{2}),(l_{1}, l_{2})} \right\}_{j_{1}, j_{2} \in \mathbb{Z}^{d}, l_{1}, l_{2} \in (\mathbb{Z}^{+})^{d}}  =  \left\{ \psi_{j_{1},l_{1}} \otimes \psi_{j_{2},l_{2}} \right\}_{j_{1}, j_{2} \in \mathbb{Z}^{d}, l_{1}, l_{2} \in (\mathbb{Z}^{+})^{d}}\] is a Wilson basis for \( L^{2}(\mathbb{R}^{2d}) \) generated by \( \psi \in M^{1,1}(\mathbb{R}) \). Set \[ k(t,y) = \sum_{j \in \mathbb{Z}^{d}} \sum_{l \in (\mathbb{Z}^{+})^{d} } \lambda_{l,j} \, \psi_{j,l}(t) \, \psi_{j,l}(y). \] Then 
\begin{align*}
& \Psi_{(n_{c(1)},n_{ c(2)}, \dots, n_{c(2d)}), (n_{c(2d+1)}, \dots, n_{c(4d)})} \\
& \indent  =  \Psi_{(n_{d+1}, \dots, n_{2d},n_{3d+1}, \dots,n_{4d}),( n_{1}, \dots, n_{d}, n_{2d+1}, \dots n_{3d})} \\
& \indent = \psi_{(n_{d+1}, \dots, n_{2d}),( n_{1}, \dots, n_{d})} \otimes  \psi_{(n_{3d+1}, \dots,n_{4d}), (n_{2d+1}, \dots n_{3d})}.
\end{align*}
By Corollary \ref{cor:isothm}
\begin{align*}
& \left\| k \right\|_{ M(c)_{v_{s}}^{p_{1},p_{2},\dots, p_{4d}}} \\
 & \asymp
\biggl( \sum_{n_{4d}\in X_{c^{-1}(4d)}} \dots \biggl( \sum_{n_{1} \in X_{c^{-1}(1)}} \left| \langle k, \Psi_{(n_{c(1)},\dots, n_{c(2d)} ),(n_{c(2d+1)}, \dots,n_{c(4d)} )} \rangle \right|^{p_{1}} \biggr)^{\frac{p_{2}}{p_{1}}} \dots \biggr)^{\frac{1}{p_{4d}}} \\
& = \biggl( \sum_{n_{4d}\in \mathbb{Z}} \biggl( \dots \biggl( \sum_{n_{2d+1} \in \mathbb{Z}^{+}} \left| \lambda_{( n_{2d+1}, \dots, n_{3d}),(n_{3d+1}, \dots, n_{4d})}  \right|^{p_{2d+1}} \biggr)^{\frac{p_{2d+2}}{p_{2d+1}}} \dots \biggr)^{\frac{p_{4d}}{p_{4d-1}}} \biggr)^{\frac{1}{p_{4d}}} \\
& = \left\| \lambda \right\|_{ \ell^{p_{2d+1}, \dots, p_{3d},p_{3d+1}, \dots p_{4d}} ((\mathbb{Z}^{+})^{d},\mathbb{Z}^{d}  )  } 
\end{align*} so \( k \in M(c)^{p_{1},p_{2},\dots, p_{4d}} \subset M(c)_{v_{s}}^{p_{1},p_{2},\dots, p_{4d}} \). The pseudodifferential operator \( A \) with kernel \( k \) has singular values equal to the elements of the sequence \( \lambda \).  Hence \( A \notin \mathcal{I}_{p}(L^{2}(\mathbb{R}^{d})) \).
\end{proof}
Notice that the proof of the previous theorem shows that Theorem \ref{thm:maintfspec} does not hold for \( p > 2 \).  That is, if \( p > 2 \) and \( k \in M(c)^{2,2, \dots, 2,p,\dots, p} \), the corresponding integral operator may not even be bounded on \( L^{2}(\mathbb{R}^{d}) \).



\section*{Acknowledgments}
The author thanks Christopher Heil for helpful discussions.

\end{document}